# The Cardinality of Infinite Games

Thomas Meyer

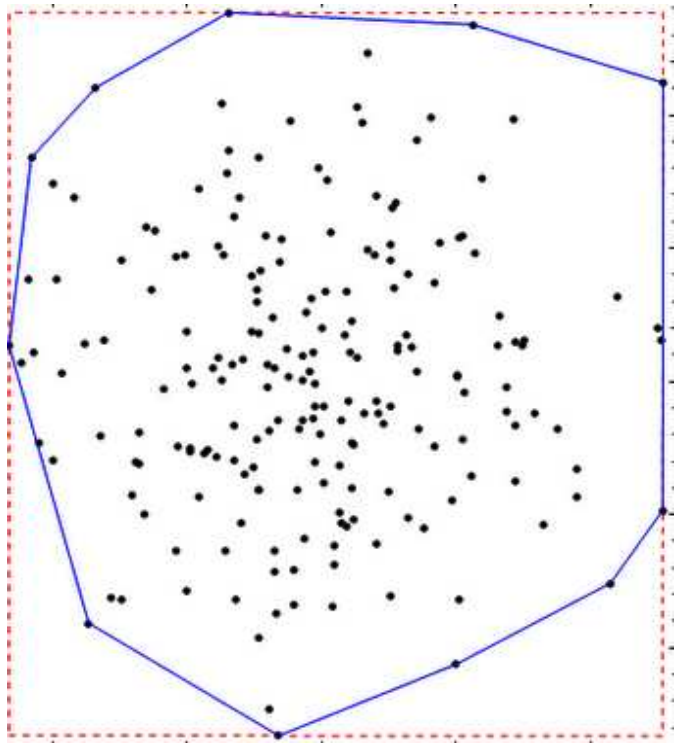

10.14.2009

Abstract

The focus of this essay is a rigorous treatment of infinite games. An infinite game is defined as a play consisting of a fixed number of players whose sequence of moves is repeated, or iterated ad infinitum. Each sequence corresponds to a single iteration of the play, where there are an infinite amount of iterations. Thus, during the first iteration of the play, the first player takes their turn, then the second, then the third, etc. Once the $n^{th}$ player's move is finished, the second iteration begins, where each player takes another turn. There are two distinct concepts within this broad definition which encompass all infinite games: the strong infinite game and the weak infinite game. Ostensibly, both satisfy the requirements to be considered an infinite game; however they differ in terms of imputations. The strong infinite game has a uniqueness qualification in that all moves must differ to the extent that no imputation (these occur at the end of any given iteration) may ever be the same. This is achieved by allowing each player a unique set of alternatives during each turn. Conversely, there is no such qualification in a weak infinite game, any payout may equal another. In point of fact, it is conceivable that in such an infinite game, all imputations may be the same (though this is a rather extreme case).

Another property shared by strong and weak infinite games (apart from their fulfilling the criterion of an infinite game) is the fact that both consist of a countably infinite amount of moves. Therefore all infinite games have a countably infinite number of turns; the set of all infinite games is composed of each strong and weak infinite game. This result is referred to as the Theorem of Infinite Games, which has a very important consequence: the ordinality of turns. That is, the moves of an infinite game have an order which they adhere to, a structure to follow. It is this structure which provides any future development or game theoretical analysis of these sorts of games with the necessary foundation.

Preface – A Brief Discussion of *Termini Technici*

      To begin, it is necessary to define the most fundamental concepts used in this essay: games, plays and turns. This paper uses the von Neumann – Morgenstern distinction between these first two constructs. The first concept is defined to be "the totality of the rules which describe it."[i] The second is stated to be "every particular instance at which the game is played."[ii] This distinction is counterintuitive, since the definition of a play given in *Theory of Games and Economic Behavior* is generally used to define the noun "game." Finally, a turn is "the occasion of a choice between various alternatives, to be made…by one of the players."[iii]

I. Introductory Remarks

      The notion that there are a finite number of turns (or moves) within an n person game is fundamental to game theory. For instance, the concept of a Nash Equilibrium is concerned with "[the definition] of equilibrium points and [proving] that a…non-cooperative game always has at least one equilibrium point."[iv] Yet, it is restricted to the case in which each of the specified *n* players has "an associated finite set of…strategies,"[v] which implies that there are solely a finite number of turns. This is consistent with the definition given by both von Neumann and Morgenstern about the nature of moves within a game. They state in the *Theory of Games and Economic Behavior*, "Consider all possible plays of the game $\Gamma$, and form the set $\Omega$ of which they are the elements…then all possible plays are simply all possible sequences [of choices between alternatives][vi] $\sigma_1,...,\sigma_v$. There exist only a finite number of such sequences, and so $\Omega$ is a finite set."[vii] Because there are a fixed number $v, v \in \mathbb{Z}^+$, of $\sigma$, there are also a finite number of turns, which they state explicitly: "We assume the number of moves – as we know that we may- to be fixed. Denote this number by *v*, and the moves themselves…by $M_1, M_2,..., M_v$."[viii] The finiteness of the number of moves has a singularly important consequence: the set of all turns is countable. From the basic principles of cardinality, there is the following theorem "A set is…countable if … it is finite."[ix] In addition, cardinality is used to define finiteness, "A set F is said to be finite if for some $n \in \mathbb{N}$ there exists a bijective mapping $\phi$ from $\{1, 2,…, n\}$ to F. The natural number *n*…is uniquely determined, and we call it the cardinality of F."[x]

      The countable nature of turns provides an inherent order and structure to the concept of a game. There exist clearly defined boundaries: each play is constituted by a specified number of turns - there is no ambiguity in this regard. In some sense, one becomes cognizant of what to expect. Yet, cardinality is not exclusively applicable to finite sets, "If X is an infinite set, there exists an injective mapping to the natural

---

[i] *Theory of Games and Economic Behavior*; Simplified Concept of a Game; Page 49.
[ii] *Theory of Games and Economic Behavior*; Simplified Concept of a Game; Page 49.
[iii] *Theory of Games and Economic Behavior*; Simplified Concept of a Game; Page 49.
[iv] *Non-Cooperative Games*; The Annals of Mathematics; Introduction; Page 286.
[v] *Non-Cooperative Games*; The Annals of Mathematics; Formal Definitions and Terminology; Page 286.
[vi] Recall that this is the definition for a turn, the "occasion" for a player to choose between alternatives.
[vii] *Theory of Games and Economic Behavior*; The Set-theoretical Description of a Game; Page 67-68.
[viii] *Theory of Games and Economic Behavior*; The Set-theoretical Description of a Game; Page 67.
[ix] *Sets, Logic and Categories*; Countable Sets; Page 24.
[x] *Mathematical Analysis: An Introduction*; Infinite Sets; Page 6.

numbers, N into X. If there exists a bijective map of ℕ into X, we say that X is countably infinite."[xi] The purpose of this paper is to treat the set of moves in a play as infinite, and to establish a bijection between these turns and ℕ. However, it may seem impractical to consider a game lasting an infinite number of terms. Additionally, as was mentioned above, most work in game theory assumes finiteness. (There is one exception, though, the infinitely repeated game. An infinite game, however, is a separate concept philosophically. It is much more concerned with the notion of an infinite number of turns within one play, whereas the purpose of infinitely repeated games is different. In these sorts of games, there are an infinite amount of plays of the given game. There are also differences in terms of particulars, which are discussed later.)

These objections are not invalid. Nevertheless, there is purpose to this paper. Proving the existence of cardinality in the case of games with an infinite number of turns (herein referred to as infinite games) serves to further strengthen the theoretical underpinnings of game theory. This will show that the inherent order obtained from countable sets still exists when the assumption of finiteness is abandoned, i.e., when a game continues ad infinitum. Furthermore, it is quite possible that this assumption may be discarded in future work in game theory, since the direction research takes is not a simple thing to predict. In this instance, research would have a secure footing given the countable nature of infinite games.

II. Formal Development of Notation[xii]

Before these considerations are proved, it is necessary to provide a rigorous definition of what exactly may be considered an infinite game. Broadly this may be defined as a play consisting of n players whose sequence of moves is repeated ad infinitum. In the notational sense, this is represented by $M_1, M_2, ..., M_n, M_{n+1}, M_{n+2}, ..., M_{2n}, M_{2n+1}, ...$. The subscript of each $M$ is of the form $an + b$, where $a, b \in \mathbb{Z}^+$. The variable a indicates the $a + 1^{st}$ iteration of the play. Additionally, b denotes the player whose turn it is. In essence, during the first iteration of the play, the first player takes their turn, then the second, then the third, etc. However, after the $n^{th}$ player's move is finished, the second iteration begins. The focus is again shifted to player one who makes another move, then again to player two, etc. This continues indefinitely. In addition, the payouts or imputations corresponding to a single iteration occur at that iteration's end. An example of the payout structure of infinite games is poker. After each hand is finished, the player with the winning cards collects the pot; the players do not wait until the play is finished to collect their winnings.

The structure of imputations for infinite games also emphasizes the dissimilarity to infinitely repeated games. For these repeated games, the payouts occur at the end of each play. However, in the infinite games, these imputations arise at the end of any iteration; were the imputations to be structured as in infinitely repeated games, no player would ever receive their payout.

Under the heading of infinite games, two new concepts may be created: strong and weak infinite games. Both of these are subsets of the set of all infinite games,

---

[xi] *Mathematical Analysis: An Introduction*; Infinite Sets; Page 6.
[xii] The notation of this paper is heavily reliant upon that of von Neumman and Morgenstern. A concise description of their notation is given on p. 50 in the *Theory of Games and Economic Behavior*.

denoted $\chi = \{\Gamma : M_1, M_2, ..., M_n, M_{n+1}, M_{n+2}, ..., M_{2n}, M_{2n+1}, ...\}$. The first is a strong infinite game, the set of which is $\lambda \subset \chi$; an element of $\lambda$ is $\eta$. Here, it is required that each turn is unique in the sense that no imputations are identical. Given the moves made by each player during the infinite number of iterations of the play, no two yield the same payout; each consists of choosing between a different set of alternatives such that no payouts are the same.

Let $i(1)_n, i(2)_n, ..., i(\alpha_n)_n$ represent the different alternatives available to player $n$ at the first turn[xiii], given $\{1, 2, ..., \alpha_n\}$. Then, $i(1')_n, i(2')_n, ..., i(\alpha'_n)_n$ describes the alternatives presented during the n$^{th}$ player's second turn, given $1', 2', ..., \alpha'_n \neq 1, 2, ..., \alpha_n$. In general, $i(1^k)_n, i(2^k)_n, ..., i(\alpha_n^k)_n$ is the set of alternatives available to $n$ at the $k + 1^{st}$ iteration of the play, with the corresponding set $\{1^k, 2^k, ..., \alpha_n^k\}$, such that $1^k, 2^k, ..., \alpha_n^k \neq 1^{k-1}, 2^{k-1}, ..., \alpha_n^{k-1}$. In addition the set $\{\sigma_1^{k-1}, ..., \sigma_n^{k-1}\}$ is defined, where $\sigma_n^k \in \{1^k, 2^k, ..., \alpha_n^k\}$, and represents the choice among alternatives made by player $n$ at the $k^{th}$ turn. Thus, it is stated that $\sigma_1^{k-1}, ..., \sigma_n^{k-1} \neq \sigma_1^k, ..., \sigma_n^k$. Finally, the payoff function for $n$ at the $k^{th}$ turn is defined as $f(\{\sigma_1^k, ..., \sigma_n^k\})_n^k$. Because no imputations may ever be the same, an additional condition must be imposed upon the set $\{\sigma_1^k, ..., \sigma_n^k\}$, and by extension $i(1^k)_n, i(2^k)_n, ..., i(\alpha_n^k)_n$, and $\{1^k, 2^k, ..., \alpha_n^k\}$. Aside from simply being unique, each $\sigma_n^k$ must be defined so that this set, when inputted into an arbitrary players payoff function, during any given iteration, outputs a unique imputation.

Conversely, in a weak infinite game, it is not necessary for each move to provide a different payout. The sole requirement is that there are an infinite number of moves for each player. Accordingly, $f(\{\sigma_1^k, ..., \sigma_n^k\})_n^k$ may equal any (or all) other imputations. This implies that $\sigma_n^k$ may be equivalent to an element(s) $\sigma_j^i$ of another set. Additionally, any one (or all) of the different alternatives for an arbitrary player $i(1)_n, i(2)_n, ..., i(\alpha_n)_n$ may be equivalent to any other possible alternative. Therefore, an element of $\{1^k, 2^k, ..., \alpha_n^k\}$, may be identical to another element from a previous or future set. In point of fact, it is possible that the given set $\{1^k, 2^k, ..., \alpha_n^k\} = \{1^{k-1}, 2^{k-1}, ..., \alpha_n^{k-1}\}$. The same can also be said for $i(1)_n, i(2)_n, ..., i(\alpha_n)_n$, or $\{\sigma_1^k, ..., \sigma_n^k\}$. However, these are the most extreme (and unlikely) outcomes. Lastly, the set of all weak infinite games is denoted $\rho \subset \chi$, an element of this set is designated by $\gamma$.

III. Statement of the Theorem[xiv]

---

[xiii] It may seem somewhat strange that there are $n$ alternatives available to each player. This is simply an extension of the von Neumann - Morgenstern notation, which specifies this number of alternatives should exist for any given player in any sort of finite game.

[xiv] Note that section III follows a more standard theorem - proof organization. This section, as distinct from the rest of the essay, invokes a more informal style, often used in proofs.

**The Theorem of Infinite Games**: Every element of the set $\chi$ is countably infinite. To be specific, $|\eta| = \aleph_0$ and $|\gamma| = \aleph_0$.

Proof:

Let us first consider the claim that $|\eta| = |N|$, for $\eta \in \lambda$. We state that $\exists$ a function $f : \mathbb{N} \to M_n \ni M_n = M_{f(\mathbb{N})}$ where $n = \mathbb{N}$. This $\vDash$, $M_{f(1)} = M_1$, $M_{f(2)} = M_2$, etc. Because $\exists$ an indefinite number of the $M_i$ with subscripts of the form $M_{an+b}$, the function $f$ is onto; $\forall x \in \mathbb{N}$, $f(x) : x \to M_x$. However, $\because \exists! M_i$ which maps to the arbitrary payout $f(\{\sigma_1^k, ..., \sigma_n^k\})_n^k \Rightarrow \exists! x$ which corresponds to this payout. Therefore, $f$ is also one to one, meaning $f$ is a bijection from $N$ to the $M_i$. By definition, $|\eta| = |N|$, or $|\eta| = \aleph_0$.

Now we turn to the claim that $|\gamma| \leq |N|$, where $\gamma \in \rho$. Again, $\exists$ a function $g : \mathbb{N} \to M_n \ni M_n = M_{f(\mathbb{N})}$ where $n = \mathbb{N}$. Using the same argument as the one above, $g$ is onto. However, we may not state that $g$ is one to one. By the definition of a weak infinite game, $\exists$ multiple $M_i$ mapping to the imputation $f(\{\sigma_1^k, ..., \sigma_n^k\})_n^k$. Consequently, we $\vdash \exists$ multiple $x$ which map to this payout through $g$.[xv] This $\vDash g : N \to M_n$ is an injection from the $M_i$ into $N \Rightarrow |\gamma| \leq |N| \Rightarrow |\gamma| = |N| = \aleph_0$.

∎

IV. Concluding Remarks

The countability of an infinite game has a singularly important consequence. Namely, the components of the game, the turns, have a well defined order. One is first, another second, etc. This order of moves exists because of the inherent link between cardinality and ordinality brought on by the current paradigm in ZFC[xvi] set theory; i.e., to begin counting the elements of a set entails first being able to order them into a sequence, an intuitive notion. This is described aptly in the following passage: "In current set theory, in order to assign a size to a set, i.e. a cardinal number, we must first be able to enumerate its elements along a well-ordered sequence, i.e., assign to it an ordinal number. This is because in ZF[C] sizes are sets (initial ordinals) which carry an inherent well-ordering…[it is known that] the theory of cardinal numbers is interwoven with that of ordinals."[xvii] This inherent ordinality implies that there is a structure contained within all infinite games, regardless of their having an indeterminate number of turns.

These references to set theory are quite appropriate since this field and the Theorem of Infinite Games accomplish the same end, though in different contexts. Set

---

[xv] This is the most average (and likely) contingency. The extreme would be plays in which all turns correspond to the same utility, implying that all $x$ did as well.

[xvi] This denotes the Zermelo – Fraenkel axioms, in addition to the axiom of choice.

[xvii] *Cardinality without Enumeration*; Introduction; Page 121.

Theory creates a foundation for all mathematics to be built upon; it describes the mechanics of the simplest mathematical objects: sets. Likewise, the cardinality of infinite games provides their turns with order, which establishes the most basic underpinning for any further game theoretical analysis. If an infinite game had no structure - no rules dictating when each player may take their turn, then all subsequent analysis or development would be well nigh impossible and impractical. No concept constructed with the finite game would be readily applicable.

        Another fitting analogy is the von Neumann axioms governing the ordinality, combining, and algebra of utility. These were necessary to establish the most basic tenants of utility, such as linearity, probabilistic weights, etc. Without such relationships firmly established, the analysis of finite games could not have continued. As is stated in the *Theory of Games and Economic Behavior*: "this possibility – i.e. the completeness of the system of (individual) preferences – must be assumed…for the purposes of the [theory of games]." [xviii] This is precisely the case for the examination or further development of infinite games: it is imperative that the Theorem of Infinite Games be used.

---

[xviii] *Theory of Games and Economic Behavior*; The Notion of Utility; Page 29.